%% file: metrics-arxiv.tex
\documentclass[11pt]{amsart}
\usepackage{amssymb,harvard}

\newtheorem{theorem}{Theorem}

\theoremstyle{definition}

\newcommand{\ignore}[1]{}

\newcommand{\munu}{(\mu,\nu)}
\newcommand{\R}{\mbox{I\hspace{-.2em}R}}
\newcommand{\Z}{{\bf Z}}
\newcommand{\id}{{\bf 1}}
\newcommand{\diam}{\mbox{diam}}

\begin{document}

\title{On choosing and bounding probability metrics}
\author{Alison L. Gibbs}
\address{Department of Mathematics and Statistics\\York University\\
4700 Keele Street\\ Toronto, Ontario, Canada, M3J 1P3}
\email{agibbs@mathstat.yorku.ca}

\author{Francis Edward Su}
\address{Department of Mathematics\\
Harvey Mudd College\\ Claremont, CA, U.S.A. 91711}
\email{su@math.hmc.edu}

\dedicatory{Manuscript version January 2002}

\thanks{First author supported in part by an NSERC postdoctoral fellowship.
Second author acknowledges the hospitality of the Cornell School
of Operations Research during a sabbatical in which this was
completed.  The authors thank Jeff Rosenthal and Persi Diaconis
for their encouragement of this project and Neal Madras for helpful
discussions.}



\begin{abstract}
When studying convergence of measures, an important issue is the
choice of probability metric.  We provide a summary and some new
results concerning bounds among some important probability
metrics/distances that are used by statisticians and
probabilists.  
Knowledge of other metrics can provide a means of
deriving bounds for another one in an applied problem.
Considering other metrics can also provide alternate insights.
We also give examples that show that rates of convergence can
strongly depend on the metric chosen.  Careful consideration is
necessary when choosing a metric.

\medskip

\noindent
{\sc Abr\'{e}g\'{e}.}
Le choix de m\'{e}trique de probabilit\'{e} est une d\'{e}cision
tr\`{e}s importante lorsqu'on \'{e}tudie la convergence des mesures.
Nous vous fournissons avec un sommaire de plusieurs m\'{e}triques/distances
de probabilit\'{e} couramment utilis\'{e}es par des statisticiens(nes)
at par des probabilistes, ainsi que certains nouveaux r\'{e}sultats qui se
rapportent \`{a} leurs bornes.  Avoir connaissance d'autres m\'{e}triques
peut vous fournir avec un moyen de d\'{e}river des bornes pour
une autre m\'{e}trique dans un probl\`{e}me appliqu\'{e}.  Le fait de
prendre en consid\'{e}ration plusieurs m\'{e}triques vous permettra
d'approcher des probl\`{e}mes d'une mani\`{e}re diff\'{e}rente.  Ainsi, nous vous
d\'{e}montrons que les taux de convergence peuvent d\'{e}pendre de
fa\c{c}on importante sur votre choix de m\'{e}trique.
Il est donc important de tout consid\'{e}rer lorsqu'on doit choisir
une m\'{e}trique.

\end{abstract}

\maketitle


\section{Introduction}
Determining whether a sequence of probability measures converges
is a common task for a statistician or probabilist.  In many
applications it is important also to \emph{quantify} that
convergence in terms of some \emph{probability metric}; hard
numbers can then be interpreted by the metric's meaning, and one
can proceed to ask qualitative questions about the nature of that
convergence.

There are a host of metrics available to quantify the distance
between probability measures; some are not even metrics in the
strict sense of the word, but are simply notions of ``distance''
that have proven useful to consider. How does one choose among all
these metrics?  Issues that can affect a metric's desirability
include whether it has an interpretation applicable to the
problem at hand, important theoretical properties, or useful
bounding techniques.

Moreover, even after a metric is chosen, it can still be useful to
familiarize oneself with other metrics, especially if one also
considers the relationships among them.  One reason is that
bounding techniques for one metric can be exploited to yield
bounds for the desired metric.  Alternatively, analysis of a
problem using several different metrics can provide complementary
insights.

The purpose of this paper is to review some of the most important
metrics on probability measures and the relationships among them.
This project arose out of the authors' frustrations in
discovering that, while encyclopedic accounts of probability
metrics are available (e.g., \citeasnoun{rachev}), relationships
among such metrics are mostly scattered through the literature or
unavailable.
Hence in this review we collect in one place
descriptions of, and bounds between, ten important probability
metrics.  We focus on these metrics because they are either
well-known, commonly used, or admit practical bounding
techniques. We limit ourselves to metrics between probability
measures ({\em simple metrics}) rather than the broader context
of metrics between random variables ({\em compound metrics}).

We summarize these relationships in a handy reference diagram
(Figure \ref{fig:metrics}), and provide new bounds between several
metrics. We also give examples to illustrate that, depending on
the choice of metric, rates of convergence can differ both
quantitatively and qualitatively.

\input{diagram}

This paper is organized as follows.  Section~\ref{metrics} reviews
properties of our ten chosen metrics. Section~\ref{relationships}
contains references or proofs of the bounds in Figure
\ref{fig:metrics}. Some examples of their applications are
described in Section~\ref{applications}. In
Section~\ref{differentrates} we give examples to show that the
choice of metric can strongly affect both the rate and nature of
the convergence.

\section{Ten metrics on probability measures}
\label{metrics}

Throughout this paper,
let $\Omega$ denote a measurable space with $\sigma$-algebra
$\mathcal{B}$. Let $\mathcal{M}$ be the space of all probability
measures on $(\Omega, \mathcal{B})$.  We consider convergence in
$\mathcal{M}$ under various notions of distance, whose
definitions are reviewed in this section.  Some of these are not
strictly metrics, but are non-negative notions of ``distance''
between probability distributions on $\Omega$ that have proven
useful in practice.
These distances are reviewed in the order given by Table
\ref{table:metrics}.

In what follows, let $\mu$, $\nu$ denote two probability measures on
$\Omega$. Let $f$ and $g$ denote their corresponding density 
functions
with respect to a $\sigma$-finite dominating measure $\lambda$
(for example, $\lambda$ can be taken to be $(\mu + \nu)/2$).
If $\Omega=\R$, let $F$, $G$ denote their corresponding distribution
functions. When needed, $X, Y$ will denote random variables on
$\Omega$ such that $\mathcal{L}(X)=\mu$ and $\mathcal{L}(Y)=\nu$.
If $\Omega$ is a metric space, it will be understood
to be a measurable space with the Borel $\sigma$-algebra.
If $\Omega$ is a bounded metric space with metric $d$, let 
$\diam(\Omega)=\sup \{ d(x,y): x,y \in \Omega \}$ denote the 
{\em diameter} of $\Omega$.

\subsection*{Discrepancy metric}
\begin{enumerate}
\item
State space: $\Omega$ any metric space.
\item
Definition:
$$
d_D(\mu,\nu) := \sup_{\mbox{\tiny all closed balls }B} | \mu(B) -
\nu(B) | .
$$
It assumes values in $[0,1]$.
\item
The discrepancy metric recognizes the metric topology of the underlying
space $\Omega$.  However, the discrepancy is scale-invariant: 
multiplying the
metric of $\Omega$ by a positive constant does not affect the
discrepancy.
\item
The discrepancy metric admits Fourier bounds, which makes it useful 
to study convergence of random walks on groups \cite[p.~34]{diaconis}.
\end{enumerate}

\subsection*{Hellinger distance}
\begin{enumerate}
\item
State space: $\Omega$ any measurable space.
\item
Definition: 
if $f$, $g$ are densities of the measures $\mu, \nu$ 
with respect to a dominating measure $\lambda$,
$$
d_H(\mu,\nu) := \left[  
\int_\Omega (\sqrt{f}-\sqrt{g})^2 \ d\lambda 
\right]^{1/2}
=
\left[  
2 \left ( 1 - \int_\Omega \sqrt{f g} \ d\lambda \right) 
\right]^{1/2}.
$$
This definition is independent of the choice of dominating measure $\lambda$.  
For a countable state space $\Omega$,
$$
d_H(\mu,\nu) := 
\left[ 
\sum_{\omega\in\Omega} 
  \left( \sqrt{\mu(\omega)} - \sqrt{\nu(\omega)} \right )^2
\right]^{1/2}
$$
\cite{diaconis-zabell}.

\item
It assumes values in $[0,\sqrt{2}]$. Some texts, 
e.g., \citeasnoun{lecam-english}, introduce a factor of 
a square root of two 
in the definition of the Hellinger distance to normalize its range of
possible values to $[0,1]$. 
We follow \citeasnoun{zolotarev-article}.
Other sources, e.g., \citeasnoun{borovkov}, \citeasnoun{diaconis-zabell}, 
define the Hellinger distance
to be the square of $d_H$.
(While $d_H$ is a metric, $d_H^2$ is not.)
An important property is that
for product measures $\mu=\mu_1 \times \mu_2$, $\nu=\nu_1 \times \nu_2$ on 
a product space $\Omega_1 \times \Omega_2$, 
$$
1- \frac{1}{2} d_H^2(\mu,\nu)
=
\left( 1- \frac{1}{2} d_H^2(\mu_1, \nu_1) \right)
\left( 1- \frac{1}{2} d_H^2(\mu_2, \nu_2) \right)
$$
\cite[p.~279]{zolotarev-article}.  
Thus one can express the distance between distributions of 
vectors with independent components in terms of the component-wise distances.  
A consequence \cite[p.~100]{reiss} of the above formula is
$d_H^2(\mu,\nu) \leq d_H^2(\mu_1,\nu_1) + d_H^2(\mu_2,\nu_2)$.

\item
The quantity $(1- \frac{1}{2} d_H^2)$ is called the {\em Hellinger affinity}.  
Apparently 
\citeasnoun{hellinger} used a similar quantity in operator theory, but
\citeasnoun[p.~216]{kakutani}
appears responsible for popularizing the Hellinger affinity and 
the form $d_H$ in his investigation of infinite products of measures.
\citeasnoun{lecam-yang} and \citeasnoun{liese-vajda}
contain further historical references.
\end{enumerate}

\subsection*{Relative entropy (or Kullback-Leibler divergence)}
\begin{enumerate}
\item
State space: $\Omega$ any measurable space.
\item
Definition: 
if $f$, $g$ are densities of the measures $\mu, \nu$ 
with respect to a dominating measure $\lambda$,
$$
d_I(\mu,\nu) := \int_{S(\mu)} f \log(f/g) \ d\lambda,
$$
where $S(\mu)$ is the support of $\mu$ on $\Omega$.  The definition is 
independent of the choice of dominating measure $\lambda$.
For $\Omega$ a countable space,
$$d_I(\mu,\nu) := \sum_{\omega \in \Omega} \mu(\omega)
   \log \frac{\mu(\omega)}{\nu(\omega)}.$$
The usual convention, based on continuity arguments, 
is to take $0 \log \frac{0}{q} = 0$ for all real
$q$ and $p \log \frac{p}{0} = \infty$ for all real non-zero $p$.
Hence the relative entropy assumes values in $[0,\infty]$.

\item
Relative entropy is not a metric, since it is not symmetric and does not
satisfy the triangle inequality. However, it has many useful
properties, including additivity over marginals of product measures:
if $\mu=\mu_1 \times \mu_2$, $\nu=\nu_1 \times \nu_2$ on 
a product space $\Omega_1 \times \Omega_2$, 
$$
d_I(\mu,\nu) = d_I(\mu_1,\nu_1) + d_I(\mu_2,\nu_2)
$$
\cite[p.~100]{cover-thomas,reiss}.

\item
Relative entropy was first defined by \citeasnoun{kullback-liebler} as
a generalization of the entropy notion of \citeasnoun{shannon}.
A standard reference on its properties is \citeasnoun{cover-thomas}.
\end{enumerate}

\subsection*{Kolmogorov (or Uniform) metric}
\begin{enumerate}
\item
State space: $\Omega=\R$.
\item
Definition:
$$
d_K(F,G) :=\sup_x |F(x)-G(x)|, \;\; x \in \R .
$$
(Since $\mu,\nu$ are measures on $\R$, it 
is customary to express the Kolmogorov metric 
as a distance between their distribution functions $F,G$.)
\item
It assumes values in $[0,1]$, and is invariant 
under all increasing one-to-one transformations of the line.
\item
This metric, due to \citeasnoun{kolmogorov-ital}, is also called the
uniform metric \cite{zolotarev-article}.
%
%
\end{enumerate}

\subsection*{L\'evy metric}
\begin{enumerate}
\item
State space: $\Omega=\R$.
\item
Definition:
$$
d_L(F,G) := \inf \{ \epsilon > 0 : G(x-\epsilon ) - \epsilon \leq
F(x) \leq G(x+\epsilon)+\epsilon , \forall x \in \R \}.
$$ 
(Since $\mu,\nu$ are measures on $\R$, it 
is customary to express the L\'evy metric 
as a distance between their distribution functions $F,G$.)
\item
It assumes values in $[0,1]$.  
While not easy to compute, the L\'evy metric does metrize weak
convergence of measures on $\R$ \cite[p.~71]{lukacs}.  
It is shift invariant, but not scale invariant.
\item
This metric was introduced by \citeasnoun[p.~199-200]{levy}.
\end{enumerate}

\subsection*{Prokhorov (or L\'evy-Prokhorov) metric}
\begin{enumerate}
\item
State space: $\Omega$ any metric space.
\item
Definition:
$$
d_P(\mu,\nu) := \inf
  \{ \epsilon > 0: \mu(B) \leq \nu(B^\epsilon)+\epsilon
\mbox{ for all Borel sets } B  \}
$$
where $ B^\epsilon = \{ x: \inf_{y \in B} d(x,y) \leq \epsilon \}
$. It assumes values in $[0,1]$.
\item
It is possible to show that this metric is symmetric in $\mu,\nu$.
See \cite[p.~27]{huber}.
\item
This metric was defined by \citeasnoun{prokhorov} 
as the analogue of the L\'evy metric for more general spaces.
While not easy to compute, this metric is theoretically important
because it metrizes weak convergence on any separable metric space
\cite[p.~28]{huber}.
Moreover, $d_P(\mu,\nu)$ is precisely the minimum distance 
``in probability'' between random variables 
distributed according to $\mu,\nu$.  This was shown by
\citeasnoun{strassen} for complete separable metric spaces and
extended by \citeasnoun{dudley1968} to arbitrary separable metric spaces. 
\end{enumerate}

\subsection*{Separation distance}
\begin{enumerate}
\item
State space: $\Omega$ a countable space.
\item
Definition:
$$
d_S(\mu,\nu) := \max_i \left( 1 - \frac{\mu(i)}{\nu(i)} \right).
$$
\item
It assumes values in $[0,1]$.  However, it not a metric.
\item
The separation distance was advocated by \citeasnoun{aldous-diaconis}
to study Markov chains because 
it admits a useful characterization in terms of strong uniform times.
\end{enumerate}

\subsection*{Total variation distance}
\begin{enumerate}
\item
State space: $\Omega$ any measurable space.
\item
Definition:
\begin{eqnarray}
\label{eq:tv-sup}
d_{TV}(\mu,\nu)
& := & \sup_{A \subset \Omega} |\ \mu (A)-\nu(A)\ | \\
\label{eq:tv-bayes}
& = & \frac{1}{2} \max_{|h| \leq 1}
   \left | \int h\ d \mu - \int h\ d \nu \right |
\end{eqnarray}
where $h : \Omega \rightarrow \R$ satisfies $|h(x)| \leq 1$. This
metric assumes values in $[0,1]$.

For a countable state space $\Omega$, the definition above becomes
$$
d_{TV} := \frac{1}{2} \sum_{x \in \Omega} |\,\mu(x)-\nu(x)\,|
$$
which is half the $L^1$-norm between the two measures. Some
authors (for example, \citeasnoun{tierney-MCMCinpractice}) define
total variation distance as twice 
this definition.


\item Total variation distance
has a coupling characterization:
$$
d_{TV}(\mu,\nu) = \inf \{ \Pr (X \neq Y): \mbox{ r.v. } X,\; Y
\mbox{ s.t. } {\mathcal{L}}(X)=\mu,\; {\mathcal{L}}(Y)=\nu \}
$$
\cite[p.~19]{lindvall}.

\end{enumerate}

\subsection*{Wasserstein (or Kantorovich) metric}
\begin{enumerate}
\item
State space:  $\R$ or any metric space.

\item
Definition:
For $\Omega=\R$, if $F,G$ are the distribution functions of
$\mu,\nu$ respectively, the 
{\em Kantorovich metric} is defined by
\begin{eqnarray*}
d_W(\mu,\nu) & :=  & \int_{-\infty}^{\infty}  |F(x)-G(x)|\,dx \\
& = & \int_0^1 | F^{-1}(t)-G^{-1}(t)  |\,dt  .
\end{eqnarray*}
Here $F^{-1},G^{-1}$ are the inverse functions of the distribution
functions $F,G$.
For any separable metric space, this is equivalent to
\begin{equation}
\label{eq:kantorovich} d_W(\mu,\nu) :=  \sup \left \{ \left| \int
h\ d\mu - \int h\ d\nu \right| : \| h \|_L \leq 1 \right \},
\end{equation}
the supremum being taken over all $h$ satisfying the Lipschitz
condition $|h(x)-h(y)| \leq d(x,y)$, where $d$ is the metric on
$\Omega$.
\item
The Wasserstein metric assumes values in $[0, \diam(\Omega)]$,
where $\diam(\Omega)$ is the diameter of the metric space $(\Omega,d)$.
This metric metrizes weak convergence on spaces of bounded diameter,
as is evident from Theorem \ref{wass-prok} below.

\item
By the Kantorovich-Rubinstein theorem, the Kantorovich metric
is equal to the {\em Wasserstein metric}:
$$
d_W(\mu,\nu) = \inf_J \{ \mbox{\bf E}[d(X,Y)] :
\mathcal{L}(X)=\mu,
 \mathcal{L}(Y)=\nu \},
$$
where the infimum is taken over all joint distributions $J$ with
marginals $\mu, \nu$. See \citeasnoun[Theorem 2]{szulga}.  
\citeasnoun[p.~342]{dudley} traces some of the history of these metrics.
\end{enumerate}

\subsection*{$\chi^2$-distance}
\begin{enumerate}
\item
State space: $\Omega$ any measurable space.
\item
Definition:
if $f$, $g$ are densities of the measures $\mu, \nu$ 
with respect to a dominating measure $\lambda$, 
and $S(\mu), S(\nu)$ are their supports on $\Omega$,
$$
d_{\chi^2}(\mu,\nu) :=
    \int_{S(\mu) \cup S(\nu)} \frac{(f - g)^2}{g} \ d\lambda.
$$
This definition is independent of the choice of dominating measure 
$\lambda$.
This metric assumes values in $[0,\infty]$.
For a countable space $\Omega$ this reduces
to:
$$
d_{\chi^2}(\mu,\nu) :=  
\sum_{\omega\in S(\mu) \cup S(\nu)}
    \frac{(\mu(\omega)-\nu(\omega))^2}{\nu(\omega)} .
$$
This distance is not symmetric in $\mu$ and $\nu$; beware that the order of
the arguments varies from author to author.  We follow
\citeasnoun{csiszar} and \citeasnoun{liese-vajda} 
because of the remarks following equation
(\ref{f-divergence}) below and the natural order of arguments
suggested by inequalities (\ref{entropy-upper-bound}) 
and (\ref{dsc-entropy}).  
\citeasnoun[p.~98]{reiss} takes the
opposite convention as well as defining the $\chi^2$-distance 
as the square root of the above expression.
\item
The $\chi^2$-distance is not symmetric, and therefore not a metric.
However, like the Hellinger distance and relative entropy,
the $\chi^2$-distance between product measures can be bounded in terms
of the distances between their marginals.
See \citeasnoun[p.~100]{reiss}.
\item
The $\chi^2$-distance has origins in mathematical statistics dating
back to Pearson.  See \citeasnoun[p.~51]{liese-vajda} for some
history.
\end{enumerate}

\bigskip
We remark that several distance notions in this section are
instances of a family of distances known as
\emph{$f$-divergences} \cite{csiszar}.  For any convex
function $f$, one may define 
\begin{equation}
\label{f-divergence}
d_f(\mu,\nu) = \sum_{\omega} 
\nu(\omega) f(\frac{\mu(\omega)}{\nu(\omega)}).
\end{equation}
Then choosing $f(x)=(x-1)^2$ yields $d_{\chi^2}$, $f(x)=x \log x$
yields $d_I$, $f(x)=|x-1|/2$ yields $d_{TV}$, and
$f(x)=(\sqrt{x}-1)^2$ yields $d_H^2$.
The family of $f$-divergences are studied in detail in
\citeasnoun{liese-vajda}.

\section{Some Relationships Among Probability Metrics}
\label{relationships}
In this section we describe in detail the relationships
illustrated in Figure~\ref{fig:metrics}.  We 
give references for relationships known to
appear elsewhere, and prove several new bounds which we state as
theorems.  In choosing the order of presentation, we
loosely follow the diagram from bottom to top.  

At the end of this section, we summarize in Theorem \ref{weak-conv}
what is known about how these metrics relate to weak convergence 
of measures.

\subsection*{The Kolmogorov and L\'evy metrics on $\R$}

For probability measures $\mu, \nu$ on $\R$ with distribution
functions $F,G$,
$$
d_L(F,G) \leq d_K(F,G).
$$
See \citeasnoun[p.~34]{huber}.  
\citeasnoun[p.~43]{petrov} notes that if $G(x)$ (i.e., $\nu$) is absolutely
continuous (with respect to Lebesgue measure), then
$$
d_K(F,G) \leq \left( 1 + \sup_x |G'(x)| \right) d_L(F,G).
$$

\subsection*{The Discrepancy and Kolmogorov metrics on $\R$}

It is evident that 
for probability measures on $\R$,
\begin{equation}
\label{u-d-bound}
d_K \leq d_D \leq 2\,d_K.
\end{equation}
This follows from the regularity of Borel sets in $\R$ and expressing
closed intervals in $\R$ as difference of rays.

\subsection*{The Prokhorov and L\'evy metrics on $\R$}
For probability measures on $\R$,
$$
d_L \leq d_P.
$$
See \citeasnoun[p.~34]{huber}.

\subsection*{The Prokhorov and Discrepancy metrics}
The following theorem shows how discrepancy may be bounded by
the Prokhorov metric by finding a suitable right-continuous function $\phi$.
For bounded $\Omega$, $\phi(\epsilon)$ gives an upper bound on the
additional $\nu$-measure of the extended ball $B^\epsilon$ over
the ball $B$, where $B^\epsilon =\{ x: \inf_{y \in B} d(x,y)
\leq \epsilon \} $.
Note that this theorem also gives an upper bound for $d_K$ 
through (\ref{u-d-bound}) above.
\begin{theorem}
Let $\Omega$ be any metric space, and
let $\nu$ be any probability measure satisfying
$$
\nu(B^\epsilon) \leq \nu(B) + \phi(\epsilon)
$$
for all balls $B$ and complements of balls $B$ and some right-continuous
function $\phi$.
Then for any other probability measure $\mu$,
if $d_P(\mu,\nu) = x$,
then $d_D(\mu,\nu)\leq x + \phi(x)$.
\end{theorem}
As an example, on the circle or line, if $\nu = U $ is the uniform
distribution, then $\phi(x)=2x$ and hence
$$
d_D(\mu,U) \leq 3 \,  d_P(\mu,U).
$$

\begin{proof}
For $\mu, \nu$ as above,
$$
\mu(B)- \nu(B^x) \geq \mu(B)- \nu(B) -\phi(x).
$$
And if $d_P(\mu,\nu) = x$,
then $\mu(B)-\nu(B^{\tilde{x}}) \leq \tilde{x}$ for all
$\tilde{x} > x $ and all Borel sets $B$.
Combining with the above inequality, we see that
$$
\mu(B)- \nu(B) -\phi(\tilde{x}) \leq \tilde{x} .
$$
By taking the supremum over $B$ which are balls or complements of
balls, obtain
$$
\sup_B (\mu(B)-\nu(B)) \leq \tilde{x} + \phi(\tilde{x}).
$$
The same result may be obtained for $\nu(B)-\mu(B)$ by noting that
$\nu(B)-\mu(B)= \mu(B^c) -\nu(B^c)$ which, after taking the supremum
over $B$ which are balls or complements of balls, obtain
$$
\sup_B (\nu(B)-\mu(B)) = \sup_{B^c} (\mu(B^c)-\nu(B^c)) \leq \tilde{x} + \phi(\tilde{x})
$$
as before.  Since the supremum over balls and complements of balls
will
be larger than the supremum over balls, if $d_P(\mu,\nu) = x$,
then $d_D(\mu,\nu)\leq \tilde{x} + \phi(\tilde{x})$ for all $\tilde{x}>x$.
For right-continuous $\phi$, the theorem follows by taking the limit
as $\tilde{x}$ decreases to $x$.
\end{proof}

\subsection*{The Prokhorov and Wasserstein metrics}
\citeasnoun[p.~33]{huber} shows that
$$
\left( d_P \right)^2 \leq d_W \leq 2\ d_P$$ for probability
measures on a complete separable metric space whose metric $d$ is
bounded by $1$.  
More generally, we show the following:
\begin{theorem}
\label{wass-prok}
The Wasserstein and Prokhorov metrics satisfy
$$ 
\left( d_P \right)^2 \leq d_W \leq (\mbox{\rm{diam}}(\Omega) +1)\ d_P
$$
where $\mbox{\rm{diam}}(\Omega)= \sup \{ d(x,y): x,y \in \Omega \}$.
\end{theorem}

\begin{proof}
For any joint distribution $J$ on random variables $X,Y$,
\begin{eqnarray*}
\mbox{\bf E}_J [d(X,Y)] 
& \leq & \epsilon \cdot \Pr (d(X,Y) \leq
\epsilon)
       + \diam(\Omega) \cdot \Pr (d(X,Y) > \epsilon) \\
& = & \epsilon + (\diam(\Omega)
    - \epsilon) \cdot \Pr (d(X,Y) > \epsilon)
\end{eqnarray*}
If $d_P(\mu,\nu) \leq \epsilon$, we can choose a coupling so that
$\Pr (d(X,Y) > \epsilon)$ is bounded by $\epsilon$
\cite[p.~27]{huber}.  Thus
$$
\mbox{\bf E}_J [d(X,Y)] 
\leq \epsilon +
  (\diam(\Omega) -\epsilon) \epsilon \leq
    (\diam(\Omega) + 1) \epsilon.
$$
Taking the infimum of both sides over all couplings, we obtain
$$
d_W \leq (\diam(\Omega)+1) d_P.
$$

To bound Prokhorov by Wasserstein, use Markov's inequality and
choose $\epsilon$ such that $d_W(\mu,\nu)=\epsilon^2$.  Then
$$
\Pr (d(X,Y)>\epsilon) \leq \frac{1}{\epsilon} \mbox{\bf
E}_J[d(X,Y)] \leq \epsilon
$$
where $J$ is any joint distribution on $X, Y$. By Strassen's
theorem \cite[Theorem 3.7]{huber},
$\Pr(d(X,Y)>\epsilon) \leq \epsilon$ is equivalent to $\mu(B)
\leq \nu(B^\epsilon)+\epsilon$ for all Borel sets $B$, giving
$d_P^2 \leq d_W$.
\end{proof}

No such upper bound on $d_W$ holds if $\Omega$ is not bounded.
\citeasnoun[p.~330]{dudley} cites the following example on $\R$.  
Let $\delta_x$ denote the delta measure at $x$.  The measures
$P_n:=((n-1)\delta_0 + \delta_n)/n$ converge to $P:=\delta_0$ under the
Prokhorov metric, but $d_W(P_n, P)=1$ for all $n$.  Thus
Wasserstein metric metrizes weak convergence only on state spaces of
bounded diameter.

\subsection*{The Wasserstein and Discrepancy metrics}
The following bound can be recovered using the bounds through
total variation (and is therefore not included on Figure
\ref{fig:metrics}), but we include this direct proof for
completeness.

\begin{theorem}
If $\Omega$ is finite,
$$
d_{\min} \cdot d_D \leq d_W
$$
where $d_{\min} = \min_{x \neq y} d(x,y)$ over distinct pairs of
points in $\Omega$.
\end{theorem}

\begin{proof}
In the equivalent form of the Wasserstein metric, Equation
(\ref{eq:kantorovich}), take
$$
h(x) = \left \{ \begin{array}{l l }
                  d_{\min}    & \mbox{for $x$ in $B$} \\
                  0           & \mbox{otherwise} \\
                \end{array}
  \right .
$$
for $B$ any closed ball.
$h(x)$ satisfies the Lipschitz condition. Then
\begin{eqnarray*}
d_{\min} \cdot  | \mu(B) - \nu(B)| & = &
\left | \int_\Omega h \, d\mu - \int_\Omega h \, d\nu \right | \\
& \leq & d_W(\mu, \nu)
\end{eqnarray*}
and taking $B$ to be the ball that maximizes $|\mu(B)-\nu(B)|$
gives the result.
\end{proof}

On continuous spaces, it is possible for $d_W$ to converge to 0 while
$d_D$ remains at 1. For example, take delta measures
$\delta_\epsilon$ converging on $\delta_0$.

\subsection*{The Total Variation and Discrepancy metrics}
It is clear that
\begin{equation}
\label{TV-D}
d_D \leq d_{TV}
\end{equation}
since total variation is the supremum over a larger class of sets
than discrepancy.

No expression of the reverse type can hold since the total
variation distance between a discrete and continuous distribution
is 1 while the discrepancy may be very small.  Further examples
are discussed in Section \ref{differentrates}.

\subsection*{The Total Variation and Prokhorov metrics}
\citeasnoun[p.~34]{huber} proves the following bound for
probabilities on metric spaces:
$$
d_P \leq d_{TV}.
$$

\subsection*{The Wasserstein and Total Variation metrics}
\begin{theorem}
The Wasserstein metric and the total variation
distance satisfy the following relation:
$$
d_W \leq \mbox{\rm{diam}}(\Omega) \cdot d_{TV}
$$
where $\mbox{\rm}{diam}(\Omega)= \sup \{ d(x,y): x,y \in \Omega \}$.
If $\Omega$ is a finite set, there is a bound the other way. 
If $d_{\min} = \min_{x \neq y} d(x,y)$ over distinct pairs of
points in $\Omega$,
then
\begin{equation}
\label{tv-wass}
 d_{\min} \cdot d_{TV}\leq d_{W} .
\end{equation}
\end{theorem}

Note that on an infinite set no such relation of the second type
can occur because $d_W$ may converge to 0 while $d_{TV}$ remains fixed
at 1. ($\min_{x \neq y} d(x,y)$ could be 0 on an infinite set.)

\begin{proof}
The first inequality follows from the coupling characterizations
of Wasserstein and total variation by taking the infimum of the
expected value
over all possible joint
distributions
of both sides of:
$$
d(X,Y) \leq \id_{ X \neq Y } \cdot \diam(\Omega).
$$
The reverse inequality follows similarly from:
$$
d(X,Y) \geq \id_{ X \neq Y  } \cdot \min_{a \neq b} d(a,b). 
$$
\end{proof}

\subsection*{The Hellinger and Total Variation metrics}
\begin{equation}
\label{hellinger-tv}
\frac{\left( d_H \right)^2}{2} \leq d_{TV} \leq d_H.
\end{equation}
See \citeasnoun[p.~35]{lecam-french}.
%

\subsection*{The Separation distance and Total Variation}
It is easy to show (see, e.g., \citeasnoun[p.~71]{aldous-diaconis})
that
\begin{equation}
\label{tv-sep}
d_{TV} \leq d_S.
\end{equation}
As Aldous and Diaconis note, there is no general reverse
inequality, since if $\mu$ is uniform on $\{1,2,...,n\}$ and $\nu$
is uniform on $\{1,2,...,n-1\}$ then $d_{TV}\munu = 1/n$ but
$d_S\munu = 1$.

\subsection*{Relative Entropy and Total Variation}
For countable state spaces $\Omega$,
$$ 2 \ \left( d_{TV} \right)^2 \leq d_I.$$
This inequality is due to \citeasnoun{kullback67}.  Some small
refinements are possible where the left side of the inequality is
replaced with a polynomial in $d_{TV}$ with more terms; see
\citeasnoun[p.~110-112]{mathai-rathie}.

\subsection*{Relative Entropy and the Hellinger distance}
$$ \left( d_H \right)^2 \leq d_I.$$
See \citeasnoun[p.~99]{reiss}.

\subsection*{The $\chi^2$-distance and Hellinger distance}
$$
d_H(\mu,\nu) \leq \sqrt{2} (\ d_{\chi^2}(\mu,\nu))^{1/4} .
$$
See \citeasnoun[p.~99]{reiss}, who also shows that 
if the measure $\mu$ is dominated by $\nu$, then the above inequality
can be strengthened:
$$
d_H(\mu,\nu) \leq ( d_{\chi^2}(\mu,\nu) )^{1/2}.
$$

\subsection*{The $\chi^2$-distance and Total Variation}
For a countable state space $\Omega$,
$$ d_{TV}\munu
=\frac{1}{2} \sum_{\omega\in\Omega} 
\frac{ |\mu(\omega)-\nu(\omega)| }{\sqrt{\nu(\omega)}} \sqrt{\nu(\omega)} 
\leq \frac{1}{2} \sqrt{ d_{\chi^2} \munu }
$$
where the inequality follows from Cauchy-Schwarz.
On a continuous state space, if $\mu$ is dominated by $\nu$ the
same relationship holds;
see \citeasnoun[p.~99]{reiss}.

\subsection*{The $\chi^2$-distance and Relative Entropy}

\begin{theorem}
The relative entropy $d_I$ and the $\chi^2$-distance $d_{\chi^2}$
satisfy
\begin{equation}
\label{entropy-upper-bound}
d_I (\mu,\nu)\leq \log\left[ 1 + d_{\chi^2}(\mu,\nu) \right].
\end{equation}
In particular, $d_I(\mu,\nu) \leq d_{\chi^2}(\mu,\nu)$.
\end{theorem}

\begin{proof}
Since log is a concave function, Jensen's inequality yields
$$
d_I(\mu, \nu)
 \leq \log\left(\int_\Omega (f/g)\ f \ d\lambda \right) \leq \log\left(1 +
 d_{\chi^2}(\mu,\nu)\right) \leq d_{\chi^2}(\mu,\nu),
$$
where the second inequality is obtained by noting that
$$
\int_\Omega \frac{(f-g)^2}{g} \ d\lambda = \int_\Omega \left(
\frac{f^2}{g} - 2f + g \right)\ d\lambda = \int_\Omega
\frac{f^2}{g} \ d\lambda - 1.
$$
\end{proof}
\citeasnoun[p.~710]{diaconis-saloffcoste-log}
derive the following alternate upper bound for the relative entropy in terms
of both the $\chi^2$ and total variation distances.
\begin{equation}
\label{dsc-entropy}
d_I\munu \leq d_{TV}\munu + \frac{1}{2} \, d_{\chi^2}\munu.
\end{equation}

\subsection*{Weak convergence}
In addition to using 
Figure \ref{fig:metrics}
to recall specific bounds, our diagram there can also be used to 
discern relationships between topologies on the space of
measures.  For instance, we can see from the mutual arrows between 
the total variation and Hellinger metrics that they 
generate equivalent topologies.  Other mutual arrows on the diagram
indicate similar relationships, subject to the restrictions
given on those bounds.

Moreover, since we know that the Prokhorov and L\'evy
metrics both metrize
weak convergence, we can also tell which other metrics metrize
weak convergence on which spaces, which we summarize in the following
theorem:

\begin{theorem}
\label{weak-conv}
For measures on $\R$, the L\'evy metric metrizes weak convergence.
Convergence under the discrepancy and Kolmogorov metrics 
imply weak convergence (via the L\'evy metric).  
Furthermore, these metrics metrize weak convergence 
$\mu_n\rightarrow \nu$ if the limiting metric $\nu$ is
absolutely continuous with respect to Lebesgue measure on $\R$.

For measures on a measurable space $\Omega$, 
the Prokhorov metric metrizes weak
convergence.  Convergence under the Wasserstein metric implies 
weak convergence. 

Furthermore, if $\Omega$ is bounded, the Wasserstein metric metrizes
weak convergence (via the Prokhorov metric), and convergence under any
of the following metrics implies weak convergence: total variation,
Hellinger, separation, relative entropy, and the $\chi^2$-distance.

If $\Omega$ is both bounded and finite, the total variation and Hellinger
metrics both metrize weak convergence.
\end{theorem}

This follows from chasing the diagram in Figure \ref{fig:metrics},
noting the existence of mutual
bounds of the L\'evy and Prokhorov metrics with other metrics 
(using the results surveyed in this section)
and reviewing conditions under which they apply.

\section{Some Applications of Metrics and Metric Relationships}
\label{applications}

We describe some of the applications of these metrics 
in order to give the reader a feel for how they 
have been used, and describe how some authors
have exploited metric relationships to obtain bounds for 
one metric via another.

The notion of {\em weak convergence} of measures is an
important concept in both statistics and probability.  For instance, 
when considering a statistic $T$ that is a functional of an empirical
distribution $F$, the ``robustness'' of the statistic under small
deviations of $F$ corresponds to the continuity of
$T$ with respect to the weak topology on the space of
measures.  See \citeasnoun{huber}.  
The L\'evy and Prokhorov metrics (and Wasserstein metric on a bounded
state space) provide quantitative ways of metrizing this topology.

However, other distances that do not metrize this topology can
still be useful for other reasons.
The total variation distance is one of the most commonly used
probability metrics, because it admits natural interpretations as
well as useful bounding techniques.  For instance, in
(\ref{eq:tv-sup}), if $A$ is any event, then total variation can
be interpreted as an upper bound on the difference of
probabilities that the event occurs under two measures.
In Bayesian statistics, the error in an expected loss function due to
the approximation of one measure by another is given (for bounded loss
functions) by the total variation distance through its
representation in equation (\ref{eq:tv-bayes}).

In extending theorems on the ergodic behavior of Markov chains on
discrete state spaces to general measurable spaces, the
total variation norm is used in a number of results 
\cite{orey,nummelin}.
More recently, total variation has found applications in 
bounding rates of
convergence of random walks on groups
(e.g., \citeasnoun{diaconis}, \citeasnoun{rosenthal-expository})
and Markov chain Monte Carlo algorithms
(e.g., \citeasnoun{tierney-annals}, \citeasnoun{MCMCinpractice}). 
Much of the success in obtaining rates of convergence 
in these settings is a result of the coupling
characterization of the total variation distance,
as well as Fourier bounds.

\citeasnoun{gibbs-wasserstein} considers a Markov chain Monte
Carlo algorithm which converges in total variation distance,
but for which coupling bounds are difficult to apply since the
state space is continuous and one must wait for random variables
to couple exactly.  The Wasserstein metric has a coupling
characterization that depends on the distance between two random
variables, so one may instead consider only the time required for
the random variables to couple to within $\epsilon$, a fact
exploited by Gibbs.  For a related example with a discrete
state space, Gibbs uses the bound (\ref{tv-wass}) to obtain total
variation bounds.

Like the total variation property (\ref{eq:tv-bayes}), 
the Wasserstein metric also represents the error in the expected value 
of a certain class of functions due to the approximation
of one measure by another, as in (\ref{eq:kantorovich}), which 
is of interest in applications in Bayesian statistics.
The fact that the Wasserstein metric is a minimal distance of two
random variables with fixed distributions has also led to its use in
the study of distributions with fixed marginals 
(e.g., \citeasnoun{fixedmarginals}).

Because the separation distance has a characterization in terms of
strong uniform times (like the coupling relationship for total variation), 
convergence of a Markov chain under the
separation distance may be studied by constructing a strong uniform
time for the chain and estimating the probability in the tail
of its distribution.  See \citeasnoun{aldous-diaconis} for such
examples; they also exploit inequality (\ref{tv-sep}) to 
obtain upper bounds on the total variation.

Similarly, total variation lower bounds may be obtained via 
(\ref{TV-D}) and lower bounds on the discrepancy metric.
A version of this metric is popular among number theorists
to study uniform distribution of sequences \cite{niederreiter}; 
\citeasnoun{diaconis} suggested its use to study random
walks on groups.
\citeasnoun{su-circle} uses the discrepancy metric to bound the
convergence time of a random walk on the circle generated by a
single irrational rotation. This walk converges weakly, but not
in total variation distance because its $n$-th step probability
distribution is finitely supported but its limiting measure is
continuous (in fact, uniform).  While the Prokhorov metric
metrizes weak convergence, it is not easy to bound.  On the other
hand, for this walk, discrepancy convergence implies weak
convergence when the limiting measure is uniform; and Fourier
techniques for discrepancy allow the calculation of quantitative
bounds.  The discrepancy metric can be used similarly to study random walks on
other homogeneous spaces, e.g., \citeasnoun{su-drunkard}.

Other metrics are useful because of their special
properties.  For instance, the Hellinger distance is convenient
when working with convergence of product measures because it
factors nicely in terms of the convergence of the components.
\citeasnoun{reiss} uses this fact and the relationships 
(\ref{hellinger-tv}) between the
Hellinger and total variation distances to obtain total
variation bounds. The Hellinger distance is also used in the
theory of asymptotic efficiency (e.g., see
\citeasnoun{lecam-english}) and minimum Hellinger distance
estimation (e.g., see \citeasnoun{lindsay}).
It is used throughout \citeasnoun{ibragimov} to quantify
the rate of convergence of sequences of consistent estimators to
their parameter.
\citeasnoun{kakutani}
gives a criterion (now known as the {\em Kakutani alternative}) 
using the Hellinger affinity to determine when infinite products of 
equivalent measures are equivalent; this has applications to
stochastic processes and can be used to show the consistency
of the likelihood-ratio test.  
See \citeasnoun{jacod-shiryaev}, \citeasnoun{williams} for applications.

\citeasnoun{diaconis-saloffcoste-log} use log-Sobolev techniques to bound the
$\chi^2$ convergence of Markov chains to their limiting distributions,
noting that these also give total variation and entropy bounds.
The $\chi^2$-distance bears its name because in the discrete case
it is the well-known $\chi^2$ statistic used, for example, 
in the classic goodness-of-fit test, e.g., see \citeasnoun[p.~184]{borovkov}.
Similarly, the Kolmogorov metric between a distribution function
and its empirical estimate is used as the test statistic 
in the Kolmogorov-Smirnov goodness-of-fit test, e.g., see
\citeasnoun[p.~336]{lehmann}.

Relative entropy is widely used
because it is a quantity that arises naturally, especially 
in information theory \cite{cover-thomas}.
Statistical applications
include proving central limit theorems 
\cite{linnik,barron}
and evaluating the
loss when using a maximum likelihood versus a Bayes density
estimate \cite{hartigan}.
In the testing of an empirical distribution against an alternative,
\citeasnoun[p.~256]{borovkov} gives the relationship of the
asymptotic behaviour of the Type II error
to the relative entropy between the empirical and alternative
distributions.
In Bayesian statistics, \citeasnoun[p.~75]{bernardo-smith}
suggest that relative entropy is the natural measure for the
lack of fit of an approximation of a distribution when preferences
are described by a logarithmic score function.

Up to a constant, the asymptotic
behaviour of relative entropy, the Hellinger, and the $\chi^2$-distance 
are identical when the ratio of the density functions is near 1 
\cite[p.~178]{borovkov}.
These three distances are used extensively
in parametric families of distributions to quantify the distance between
measures from the same family indexed by different parameters.
\citeasnoun[pp.~180-181]{borovkov} shows how these
distances are related to the Fisher information in the limit
as the difference in the parameters goes to zero.

\section{Rates of Convergence that Depend on the Metric}
\label{differentrates}

We now illustrate the ways in which the choice of metric can
affect rates of convergence in one context: the convergence of a
random walk to its limiting distribution. Such examples point to
the need for practitioners to choose a metric carefully when
measuring convergence, paying attention to that metric's
qualitative and quantitative features.  We give several examples
of random walks whose qualitative convergence behavior depend
strongly on the metric chosen, and suggest reasons for this
phenomenon.

As a first basic fact, it is possible for convergence to occur in
one metric but not another.
An elementary example is the convergence of a standardized
Binomial $(n,p)$ random variable with distribution $\mu_n$ which
converges to the standard normal distribution, $\nu$, as $n
\rightarrow \infty$. For all $n < \infty$, $d_{TV}(\mu_n,\nu)=1$,
while $d_D(\mu_n,\nu) \rightarrow 0$ as $n \rightarrow \infty$.
In the random walk context, \citeasnoun{su-circle} shows that a random
walk on the circle generated by an irrational rotation converges
in discrepancy, but not total variation.  The latter fact follows
because the $n$-th step probability distribution is finitely
supported, and remains total variation distance 1 away from its
continuous limiting distribution.

However, more interesting behavior can arise.  Below, we cite a
family of random walks on a product space, indexed by some parameter, 
which not only converges under each of total variation,
relative entropy, and the $\chi^2$-distance, but exhibits
\emph{different rates of convergence} as a function of the
parameter.

We then cite another example of a family of walks that not only has different
convergence rates under two different metrics, but also exhibits
\emph{qualitatively different} convergence behavior.  This family
exhibits a {\em cutoff phenomenon} under the first metric but only
exponential decay under the second.

\subsection*{Example: convergence rates that depend on the metric}

The following family of random walks show that convergence rates
in total variation, relative entropy, and $\chi^2$-distance may
differ.  
Recall that a family of random walks, indexed by some parameter
$n$, is said to \emph{converge in $f(n)$ steps} using some
metric/distance if that metric/distance can be uniformly bounded
from above after $f(n)$ steps and is uniformly bounded away from
zero before $f(n)$ steps.  

Let $G=\Z \mod g$, a finite group with $g$ elements.  Then $G^n$
is the set of all $n$-tuples of elements from $G$.  Consider the
following continuous-time random walk on $G^n$: start at
$(0,0,...,0)$ and according to a Poisson process running at rate
$1$, pick a coordinate uniformly at random and replace that
coordinate by a uniformly chosen element from $G$. (Thus each
coordinate is an independent Poisson process running at rate
$\frac{1}{n}$.) 

\citeasnoun{su-thesis} proves that if $g$ grows with $n$ exponentially
according to $g=2^n$, then 
the random walk on $G^n$ described above converges
in $2 n\log n$ steps under the relative entropy distance, and $n^2
\log 2$ steps under the $\chi^2$-distance. However, it converges
in at most $n\log n$ steps under total variation.

In this example, 
the relative entropy and $\chi^2$-distance are calculated by using
properties of these metrics, while the 
total variation upper bound in this example is, in fact,
calculated via its metric relationship (\ref{tv-sep})
with the separation distance.

Why does the difference in rates of convergence occur?  This is due to the
fact that total variation distance, unlike relative entropy or the
$\chi^2$-distance, is insensitive to big or small elementwise
differences when the total of those differences remains the same.
For instance, consider the follow measures $\mu, \nu$ on
$\Z_{10}$. Let
$$
\mu(i)= \left\{
\begin{array}{ll}
0.6 & $ if $i=0 \\
0.1 & $ if $i=1,2,3,4 \\
0  & $ else $\end{array} \right.
 \qquad
 \nu(i)=\left\{ \begin{array}{ll}
 0.2 & $ if $i=0,1,2,3,4 \\
             0  & $ else $
             \end{array}
             \right.
.
$$
Let $U$ be the uniform distribution on $\Z_{10}$.  In this example
$d_I(\mu,U)\approx 1.075$ and $d_I(\nu,U)= 0.693$, but the total
variation distances are the same:  $d_{TV}(\mu,U) = d_{TV}(\nu,U) =
0.5$. This is true because we could redistribute mass among the
points on which $\mu$ exceeded $\nu$ without affecting the total
variation distance.

Thus it is possible for an unbalanced measure (with large
elementwise differences from uniform) to have the same total
variation distance as a more balanced measure.  The random walk
on the product space has a natural feature which rewards total
variation distance but hinders relative entropy--- the
randomization in each coordinate drops the variation distance
quickly, but balancing the measure might take longer.

\subsection*{Example: qualitatively different convergence behavior}

\citeasnoun{chung-diaconis-graham} study a family of walks that
exhibits a {\em cutoff phenomenon} under total variation; the
total variation distance stays near 1 but drops off sharply and
steeply to zero after a certain cutoff point.  The family of random walks is
on the integers mod $p$ where $p$ is an odd number, with a
randomness multiplier: the process is given by $X_0=0$ and
$X_n=2X_{n-1}+\epsilon_n \; (\mbox{mod}\, p)$ where the
$\epsilon_i$ are i.i.d. taking values $0$, $\pm 1$ each with
probability $1/3$. The stationary distribution for this process
is uniform. Using Fourier analysis,
\citeasnoun{chung-diaconis-graham} show for the family of walks
when $p=2^t-1$ that $O( \log_2 p \log \log_2 p)$ steps are
sufficient and necessary for total variation convergence, for $t$
a positive integer.  There is, in fact, a cutoff phenomenon.


However, as proven in \citeasnoun[pp.~29-31]{su-thesis}, $O(\log
p)$  steps are sufficient for convergence in discrepancy.
Moreover, the convergence is qualitatively different from total
variation, because there is no cutoff in discrepancy.  Rather,
there is only exponential decay.

Again, analysis under these two metrics sheds some light on what
the walk is doing.  For $p=2^t -1$, the ``doubling'' nature of
the process keeps the walk supported very near powers of 2. The
discrepancy metric, which accounts for the topology of the space,
falls very quickly because the walk has spread itself out even
though it is supported on a small set of values. However, the
total variation does not ``see'' this spreading; it only falls
when the support of the walk is large enough.

Thus, the difference in rates of convergence sheds some light on
the nature of the random walk as well as the metrics themselves.

\bibliographystyle{dcu}
\bibliography{metrics}

\end{document}

%% file: diagram.tex
\begin{figure}

\setlength{\unitlength}{1mm}

\begin{center}

\begin{picture}(75,180)(-37.5,-110)

\thicklines


\put(-52,32){\framebox(20,20){\Huge S}}
\put(-11,11){\vector(-1,1){20}}
\put(-26.5,27.5){\scriptsize{$x$}}

\put(-10,-10){\framebox(20,20){\Huge TV}}
\put(11,1){\vector(1,0){20}}
\put(27,2.5){\scriptsize{$x$}}
\put(31,-1){\vector(-1,0){20}}
\put(12,-5){\scriptsize{$\sqrt{2x}$}}
\put(32,-10){\framebox(20,20){\Huge H}}

\put(32,-52){\framebox(20,20){\Huge W}}

\put(62,-39){on}
\put(62,-44){metric}
\put(62,-49){spaces}

\put(10,-12){\vector(1,-1){20}}
\put(17.5,-30){\scriptsize{$x/d_{\min}$}}
\put(32.5,-30.5){\vector(-1,1){20}}
\put(17.5,-13){\scriptsize{$\mbox{diam}\Omega \cdot x$}}

\put(11,-41){\vector(1,0){20}}
\put(24,-40){\scriptsize{$\sqrt{x}$}}
\put(31,-43){\vector(-1,0){20}}
\put(11,-46.5){\scriptsize{$(\mbox{diam}\Omega+1) \, x$}}

\put(11,11){\vector(1,1){20}}
\put(26,25){\scriptsize{$\sqrt{x}/2$}}
\put(25,22){\scriptsize{$\nu$ dom $\mu$}}

\put(31,11){\vector(-1,1){20}}
\put(15.5,27.5){\scriptsize{$\sqrt{x}$}}

\put(0,11){\vector(0,1){20}}
\put(-10,26){\scriptsize{$\sqrt{x/2}$}}
\put(-10,32){\framebox(20,20){\Huge I}}


\put(11,42){\vector(1,0){20}}
\put(17,43.5){\scriptsize{$\log(1+x)$}}
\put(32,32){\framebox(20,20){\Huge $\chi^2$}}

\put(62,43){non-metric}
\put(62,38){distances}

\put(42,11){\vector(0,1){20}}
\put(45,26){\scriptsize{$\sqrt{x}$}}
\put(44,23){\scriptsize{$\nu$ dom $\mu$}}


\put(0,-31){\vector(0,1){20}}
\put(1.5,-16){\scriptsize{$x$}}
\put(-10,-52){\framebox(20,20){\Huge P}}


\put(-31,-31){\vector(1,1){20}}
\put(-18,-14){\scriptsize{$x$}}

\put(-52,-52){\framebox(20,20){\Huge D}}

\put(-43,-53){\vector(0,-1){20}}
\put(-47.5,-70.5){\scriptsize{$2x$}}
\put(-41,-73){\vector(0,1){20}}
\put(-39,-57){\scriptsize{$x$}}

\put(-52,-94){\framebox(20,20){\Huge K}}

\put(-10,-94){\framebox(20,20){\Huge L}}

\put(62,-86){on $\R$}

\put(0,-73){\vector(0,1){20}}
\put(1.5,-56){\scriptsize{$x$}}


\put(-11,-86){\vector(-1,0){20}}
\put(-29,-87.5){\scriptsize{$x$}}
\put(-31,-84){\vector(1,0){20}}
\put(-30,-83){\scriptsize{$(1+\sup|G'|)x$}}


\put(-31,-42){\vector(1,0){20}}
\put(-22.5,-40.5){\scriptsize{$x+\phi(x)$}}

\end{picture}

\end{center}
\caption{Relationships among probability metrics.  A directed
arrow from $A$ to $B$ annotated by a function $h(x)$ means that $d_A \leq
h(d_B)$. The symbol
$\mbox{diam } \Omega $ 
denotes the diameter of the probability space $\Omega$; bounds
involving it are only useful if $\Omega$ is bounded.  For $\Omega$
finite, $d_{\min} = \inf_{x,y \in \Omega} d(x,y)$.  The probability
metrics take arguments $\mu, \nu$;  ``$\nu$ dom
$\mu$'' indicates that the given bound only holds if $\nu$ dominates
$\mu$.  Other notation
and restrictions on applicability
are discussed in Section \ref{relationships}.
 \label{fig:metrics} }
\end{figure}
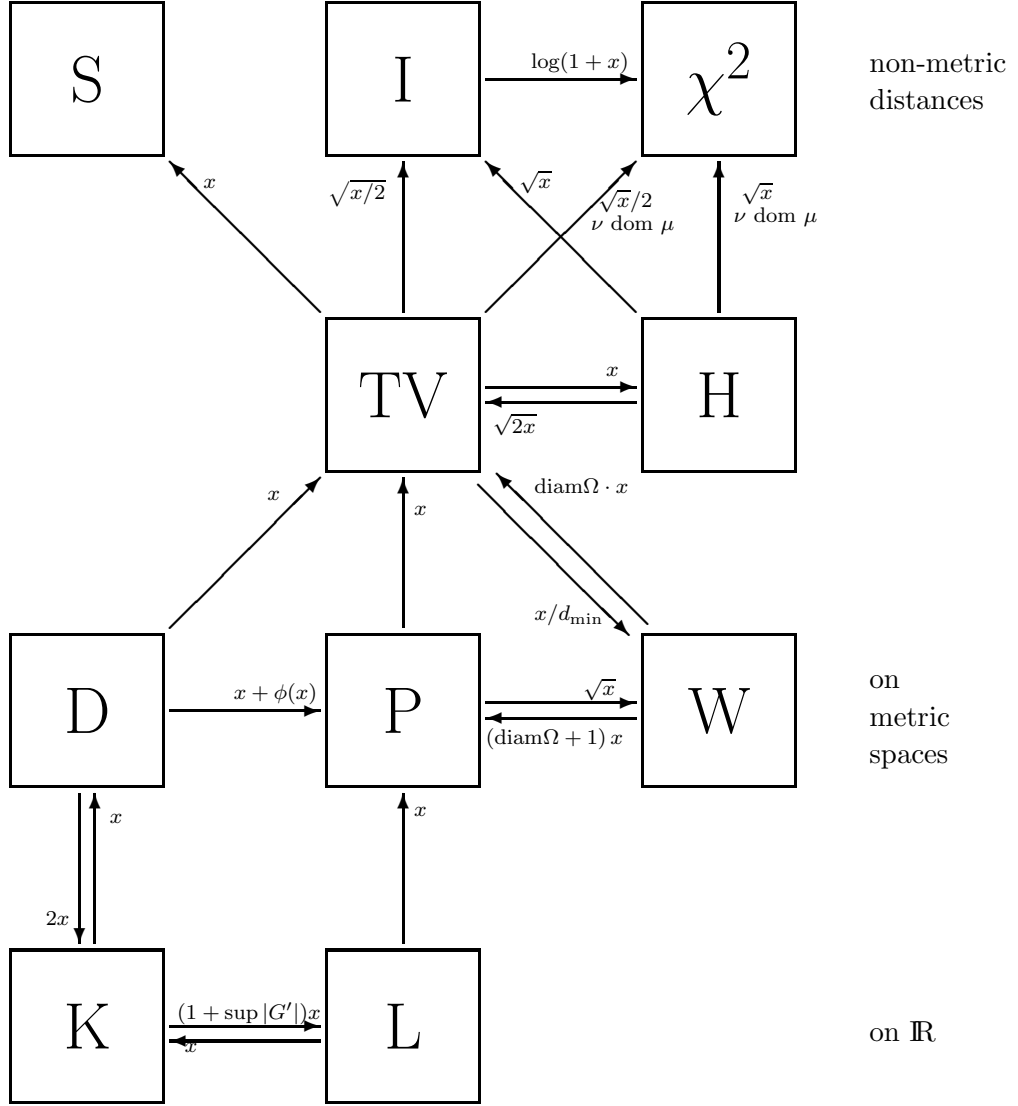

\begin{table}
\label{table:metrics}
\begin{center}
\begin{tabular}{| c | c  |}
\hline
Abbreviation & Metric \\
\hline
\hline
D & Discrepancy \\
H & Hellinger distance \\
I & Relative entropy (or Kullback-Leibler divergence) \\
K & Kolmogorov (or Uniform) metric \\
L & L\'evy metric \\
P & Prokhorov metric \\
S & Separation distance \\
TV & Total variation distance \\
W & Wasserstein (or Kantorovich) metric \\
$\chi^2$ & $\chi^2$ distance \\
\hline
\end{tabular}
\end{center}
\caption{Abbreviations for metrics used in Figure
\ref{fig:metrics}.}
\end{table}